\documentclass[12pt]{amsart}
\usepackage{amssymb}

\renewcommand{\)}{\right)}

\newcommand{\<}{\langle}
\renewcommand{\>}{\rangle}
\newcommand{\x}{\times}
\renewcommand{\bar}{\overline}
\newcommand{\abs}[1]{\left\lvert#1\right\rvert}
\newcommand{\norm}[1]{\left\lVert#1\right\rVert}
\newcommand{\st}{\:|\:}

\newcommand{\CC}{{\mathbb{C}}}
\newcommand{\RR}{{\mathbb{R}}}

\renewcommand{\phi}{\varphi}
\newcommand{\p}{\rho}

\newcommand{\tr}{{\mathrm{tr}}}

\newcommand{\B}{{\mathcal{B}}}

\renewcommand{\H}{{\mathcal{H}}}
\newcommand{\Z}{{\mathcal{Z}}}

\theoremstyle{plain}
\newtheorem{thm}{Theorem}[section]
\newtheorem{lem}[thm]{Lemma}
\newtheorem{prop}[thm]{Proposition}
\newtheorem{cor}[thm]{Corollary}

\theoremstyle{definition}
\newtheorem{defn}[thm]{Definition}

\theoremstyle{remark}

\title{{Benedicks'} Theorem for the Weyl Transform}

\author{M.~K.~Vemuri}
\address{Department of Mathematics, West Virginia University,
Morgantown, WV 26506}

\begin{document}

\begin{abstract}
If the set of points where a function is nonzero is of finite measure, and
its Weyl transform is a finite rank operator, then the function is identically
zero.
\end{abstract}

\maketitle


\section{Introduction}\label{S:intro}
That a nonzero function and its Fourier transform cannot both be sharply
localized is known as the {\em Uncertainty Principle} in Harmonic Analysis.
There are many different precise formulations of this principle, depending
on the way in which localization is quantified.  One such is Benedicks'
theorem \cite{Benedicks}: if $f \in L^1(\RR)$ and the sets
$\{x\in\RR \st f(x) \ne 0\}$ and $\{\xi\in\RR \st \hat{f}(\xi) \ne 0\}$
both have finite Lebesgue measure, then $f \equiv 0$.  In this paper, we
prove an analogue of Benedicks' theorem for the Weyl transform.  The Weyl
transform is the essence of the group Fourier transform on the Heisenberg
group, and also plays a role in the theory of pseudo-differential operators.

Let $\H=L^2(\RR)$, and $\B(\H)$ the set of bounded operators on $\H$.
If $f \in L^1(\RR^2)$, the {\em Weyl transform} of $f$ is the operator
$W(f) \in \B(\H)$ defined by
\begin{equation*}
(W(f)\phi)(t) = \iint f(x,y) e^{\pi i (xy+2yt)} \phi(t+x) \, dx dy.
\end{equation*}
Our analogue of Benedicks' theorem is the following.
\begin{thm}\label{T:main-theorem}
If the set $\{w \in \RR^2 \st f(w) \ne 0\}$ has finite Lebesgue measure and
$W(f)$ is a finite rank operator, then $f \equiv 0$.
\end{thm}

Recall that the Heisenberg group $G$ is the set of triples
\begin{equation*}
\{(x,y,z) \st x,y \in \RR, z \in \CC, \abs{z}=1\}
\end{equation*}
with multiplication defined by
\begin{equation*}
(x,y,z)(x',y',z')=\(x+x',y+y', zz'e^{\pi i(xy'-yx')}\).
\end{equation*}
According to the {\em Stone-von Neumann Theorem}, there is a unique
irreducible unitary representation $\p$ of $G$ such that
\begin{equation*}
\p(0,0,z)= zI.
\end{equation*}
The standard realization of this representation is on the Hilbert space
$\H$ by the action
\begin{equation*}
(\p(x,y,z)\phi)(t)=ze^{\pi i(xy+2yt)}\phi(t+x).
\end{equation*}
Thus, the Weyl transform may be expressed as
\begin{equation*}
W(f) = \iint f(x,y) \p(x,y,1) \, dx dy, \qquad f \in L^1(\RR^2).
\end{equation*}

If $X$ is a trace class operator on $\H$, the modified Fourier-Wigner transform
of $X$ is the function $\alpha(X):\RR^2 \to \CC$ defined by
\begin{equation*}
\alpha(X)(x,y)=\tr(X\p(x,y,1)^*).
\end{equation*}

It is well known (see e.g.\ \cite{Folland}) that if
$f \in L^1(\RR^2)$ and $W(f)$ is a trace class
operator then $\alpha(W(f))=f$, and that if $X$ is a trace class operator
on $\H$ and $\alpha(X) \in L^1(\RR^2)$ then $W(\alpha(X))=X$.  Thus we may
reformulate Theorem \ref{T:main-theorem} as

\begin{thm}\label{T:reformulation}
If $X$ is a finite rank operator on $\H$ and the set
$\{w \in \RR^2 \st \alpha(X)(w) \ne 0\}$ has finite measure, then $X=0$.
\end{thm}

Theorem \ref{T:main-theorem} was proved in \cite{NR} with the hypothesis
that $f$ is compactly supported.  However, as pointed out by Benedicks in
\cite{Benedicks},
this is a very strong hypothesis, because the support of $f$ is the
{\em closure} of the set $\{w \in \RR^2 \st f(w) \ne 0\}$, and there are open
sets of arbitrarily small measure whose closure is all of $\RR^2$.
The clever argument used in \cite{NR} doesn't seem to be
susceptible to generalization.  Our argument is closer in spirit to
Benedicks' original argument, but depends on a deep result of
Linnel \cite{Linnel}.

\section{The double induced realization}

Let $\pi:G \to \RR^2$ be the projection $\pi(x,y,z)=(x,y)$.  Then
$\pi$ is a homomorphism and $\ker(\pi)=Z(G)$, the center of $G$.
Let $U(1)=\{z \in \CC \st \abs{z}=1\}$.  For $(x,y) \in \RR^2$,
let $s(x,y)=(x,y,1)$.  Then $s$ is a section of $\pi$, i.e.
$\pi\circ s=\mathrm{id}_{\RR^2}$, and
\begin{equation*}
s(x,y) s(x',y')=(x+x',y+y',e^{\pi i (xy'-yx')})=\psi((x,y),(x',y'))s(x+x',y+y'),
\end{equation*}
where $\psi((x,y),(x',y')) \in Z(G)$ is defined by
\begin{equation*}
\psi((x,y),(x',y'))=(0,0,e^{\pi i (xy'-yx')})
\end{equation*}
Define
$e:\RR^2 \x \RR^2 \to U(1)$ by
\begin{equation*}
e((x,y), (x',y'))=e^{2\pi i (xy'-yx')}.
\end{equation*}
Note that if $g,g' \in G$ then
\begin{equation*}
gg'g^{-1}g'^{-1}=(0,0,e(\pi(g),\pi(g')).
\end{equation*}

The pairing $e$ is an alternating bicharacter, and it is perfect in the sense
that it establishes Pontryagin duality between $\RR^2$ and
$\widehat{\RR^2}$.
Let $N \subseteq \RR^2$ be a lattice (necessarily cocompact) which is
isotropic, i.e. for all $n,n' \in N$, we have $e(n,n')=1$.
Note that $N \subseteq M \subseteq N^\perp$ where
$N^\perp =\{w \in \RR^2 \st e(n,w)=1 \quad \forall n \in N\}$
and $M$ is a maximal isotropic subgroup.  Moreover,
$[N^\perp:M]=[M:N]=\mathrm{area}(\RR^2/N) = a < \infty$.
The bicharacter $e$ also establishes Pontryagin dualities
$N^\perp \cong \widehat{\RR^2/N}$, $N\cong \widehat{\RR^2/N^\perp}$, and
$M^\perp \cong \widehat{\RR^2/M}$.


There exists a character $\zeta:\pi^{-1}(M) \to U(1)$ such that
$\zeta(0,0,z)=z$ (see \cite[p3]{Theta-III}).
Let $\tau=\mathrm{Ind}_{\pi^{-1}(M)}^{\pi^{-1}(N^\perp)} \zeta$.
Then $\tau$ is an a-dimensional unitary representation of
$\pi^{-1}(N^\perp)$, and by the Stone-von Neumann-Mackey theorem
\cite[Theorem 1.2]{Theta-III},
$\mathrm{Ind}_{\pi^{-1}(N^\perp)}^G \tau \cong \mathrm{Ind}_{\pi^{-1}(M)}^G \zeta
\cong \p$, so $\tau$ is irreducible.  Note that if $n \in N$, then
$\tau(n)=\zeta(n)I$.
Let $S \subseteq N^\perp$ be a fundamental domain for the action of
$N$ on $N^\perp$, i.e. a set of coset representatives.  Then $S$ has
$a^2$ elements.  Let $\H_\tau$ be the Hilbert space on which $\tau$ acts.

\begin{lem}\label{L:pointwise-trace}
For all $T \in \B(\H_\tau)$, we have
\begin{equation*}
\sum_{n \in S} \tau(s(n)) T \tau(s(n))^{-1}=
a \, \mathrm{tr}(T) I.
\end{equation*}
\end{lem}

\begin{proof}
Let
\begin{equation*}
T'=\sum_{n \in S} \tau(s(n)) T \tau(s(n))^{-1}.
\end{equation*}
Let $g \in \pi^{-1}(N^\perp)$ and $n \in S$, then
$g=\eta s(\pi(g))$ for some $\eta \in Z(G)$, and there exist unique elements
$n' \in S$ and $m \in N$ such that $\pi(g)+n=m+n'$.  Therefore,
\begin{equation*}
\begin{aligned}
\tau(g) \tau(s(n))
=&\; \tau(g s(n))\\
=&\; \tau(\eta s(\pi(g))s(n))\\
=&\; \tau(\eta \psi(\pi(g),n)s(\pi(g)+n))\\
=&\; \tau(\eta \psi(\pi(g),n)s(m+n'))\\
=&\; \tau(\eta \psi(\pi(g),n)\psi(m,n')^{-1}s(m)s(n'))\\
=&\; \zeta(\eta \psi(\pi(g),n)\psi(m,n')^{-1}s(m))\tau(s(n'))\\
\therefore \tau(g) \tau(s(n)) T \tau(s(n))^{-1} \tau(g)^{-1}
=&\; \tau(s(n'))T\tau(s(n'))^{-1}
\end{aligned}
\end{equation*}
For $g \in \pi^{-1}(N^\perp)$ and $n \in S$, the assignment $(g,n) \mapsto n'$,
where $\pi(g)+n=m+n'$ with $n' \in S$ and $m \in N$ defines a right-action of
$\pi^{-1}(N^\perp)$ on $S$; in particular, for a fixed $g$ the map
$n \mapsto n'$ is a bijection of $S$. It follows that
$\tau(g) T' \tau(g)=T'$ for each $g \in \pi^{-1}(N^\perp)$.  It follows from
Schur's lemma that there exists $\lambda \in \CC$ such that $T'=\lambda I$.
The result follows by taking traces of both sides.
\end{proof}

Let $\B(\H_\tau)$ denote the set of all linear operators on $\H_\tau$.
If $A,B \in \B(\H_\tau)$, we put $\<A,B\>=\tr(AB^*)$, and $\norm{A}^2=\<A,A\>$.
Then $\B(\H_\tau)$ becomes a Hilbert space.

\begin{lem}
For all $n' \in N^\perp$, we have
\begin{equation*}
\tr(\tau(s(n'))=
\begin{cases}
\zeta(s(n')) a & \text{if $n' \in N$},\\
0 & \text{if $n' \notin N$}.
\end{cases}
\end{equation*}
\end{lem}

\begin{proof}
Consider $\tilde{e}: N^\perp/N \x N^\perp/N \to U(1)$ defined by
$\tilde{e}(n+N,n'+N)=e(n,n')$.  If $n' \notin N$, then
$\tilde{e}(\cdot,n'+N)$ is a non-trivial character of the finite group
$N^\perp/N$, so
\begin{equation*}
\sum_{n \in S} e(n,n') = \sum_{n \in S} \tilde{e}(n+N,n'+N)
\sum_{\bar{n} \in N^\perp/N} \tilde{e}(\bar{n}, n'+N) = 0.
\end{equation*}
If $n' \in N$, then $e(n,n') = 1$ for all $n \in N^\perp$, so
\begin{equation*}
\sum_{n \in S} e(n,n') = \abs{S}=a^2.
\end{equation*}

Let $n' \in N^\perp$ and $n \in S$.  Then
\begin{equation*}
\begin{aligned}
\tau(s(n)) \tau(s(n')) \tau(s(n))^{-1}
=&\;\tau(s(n)s(n')s(n)^{-1}s(n')^{-1}s(n'))\\
=&\;\tau((0,0,e(n,n'))s(n'))\\
=&\; e(n,n')\tau(s(n'))\\
\therefore \sum_{n \in S} \tau(s(n)) \tau(s(n')) \tau(s(n))^{-1}
=&\;
\begin{cases}
a^2 \zeta(s(n'))I & \text{if $n' \in N$},\\
0    & \text{if $n' \notin N$}.
\end{cases}
\end{aligned}
\end{equation*}
The result follows from Lemma \ref{L:pointwise-trace}
\end{proof}

\begin{cor}\label{C:finite-inner-products}
For $n, n' \in N^\perp$,
\begin{equation*}
\<\tau(s(n)), \tau(s(n'))\>=
\begin{cases}
a\zeta(\psi(n,-n')s(n-n')) & \text{if $n-n' \in N$},\\
0 & \text{if $n-n' \notin N$}.
\end{cases}
\end{equation*}
\end{cor}

\begin{proof}
\begin{equation*}
\begin{aligned}
\<\tau(s(n)), \tau(s(n'))\>
=&\; \tr(\tau(s(n))\tau(s(n'))^*)\\
=&\; \tr(\tau(s(n)s(n')^{-1}))\\
=&\; \tr(\tau(s(n)s(-n')))\\
=&\; \tr(\tau(\psi(n,-n')s(n-n')))\\
=&\; \zeta(\psi(n,-n'))\tr(\tau(s(n-n')))\\
=&\;
\begin{cases}
a\zeta(\psi(n,-n')s(n-n')) & \text{if $n-n' \in N$},\\
0 & \text{if $n-n' \notin N$}.
\end{cases}
\end{aligned}
\end{equation*}
\end{proof}

\begin{cor}\label{C:finite-onb}
The set $\{\frac{1}{\sqrt{a}}\tau(s(n)) \st n \in S\}$ is an orthonormal basis
for $\B(\H_\tau)$.
\end{cor}

\begin{proof}
Let $n,n' \in S$.  Note that $n-n' \in N$ iff $n=n'$.  Then
\begin{equation*}
\<\tau(s(n)), \tau(s(n'))\>=
\begin{cases}
a & \text{if $n = n'$},\\
0 & \text{if $n \ne n'$},
\end{cases}
\end{equation*}
so the set is orthonormal.  Since the cardinality of the set
equals the dimension of $\B(\H_\tau)$, it must be a basis.
\end{proof}

Let $\Omega_{N^\perp} \subseteq \RR^2$ be a (measurable) fundamental domain for
the action of $N^\perp$ on $\RR^2$.
Note
that $s(\Omega_{N^\perp})$ is a fundamental domain for the
action of $\pi^{-1}(N^\perp)$ on $G$.

The Hilbert space $\H_N$ for $\p_N=\mathrm{Ind}_{\pi^{-1}(N^\perp)}^G \tau$ consists of
functions $\phi:G \to \H_\tau$ such that
$\phi(hg)=\tau(h)\phi(g)$ for $h\in\pi^{-1}(N^\perp)$ and almost every
$g\in G$ with the condition that
\begin{equation*}
\norm{\phi}^2=\int_{\Omega_{N^\perp}} \norm{\phi(s(w))}^2 dw < \infty,
\end{equation*}
and $\p_N$ acts by
\begin{equation*}
(\p_N(g)\phi)(g')=\phi(g'g).
\end{equation*}
%

Let $L^0(G//\pi^{-1}(N^\perp), \B(\H_\tau))$ consist of all
measurable functions $F: G \to \B(\H_\tau)$ such that
$F(hg)=\tau(h)F(g)\tau(h)^{-1}$ for all $h \in \pi^{-1}(N^\perp)$ and
almost every $g\in G$.
For $p\in[1,\infty)$, let $L^p(G//\pi^{-1}(N^\perp), \B(\H_\tau))$ consist of
those $F \in L^0(G//\pi^{-1}(N^\perp), \B(\H_\tau))$ for which
\begin{equation*}
\norm{F}_p^p=\int_{\Omega_{N^\perp}} \norm{F(s(w))}^p \, dw < \infty.
\end{equation*}
Let $L^\infty(G//\pi^{-1}(N^\perp), \B(\H_\tau))$ consist of
those $F \in L^0(G//\pi^{-1}(N^\perp), \B(\H_\tau))$ for which
\begin{equation*}
\norm{F}_\infty=\mathrm{ess\, sup}_{g \in G} \norm{F(g)} < \infty.
\end{equation*}

Then
$L^\infty(G//\pi^{-1}(N^\perp), \B(\H_\tau)) \subseteq
L^2(G//\pi^{-1}(N^\perp), \B(\H_\tau)) \subseteq
L^1(G//\pi^{-1}(N^\perp), \B(\H_\tau))$.
Moreover, $L^2(G//\pi^{-1}(N^\perp), \B(\H_\tau))$ is a Hilbert space with
the inner product
$$
\<F,F'\>=\int_{\Omega_{N^\perp}} \<F(s(w)),F'(s(w))\> \, dw
$$

For $n \in N^\perp$, define $\Xi_n: G \to \B(\H_\tau)$ by
\begin{equation*}
\Xi_n(g')=e(\pi(g'),n)\tau(s(n)).
\end{equation*}
Then $\Xi_n \in L^\infty(G//\pi^{-1}(N^\perp), \B(\H_\tau))$.

\begin{prop}\label{P:trig-ops}
For $n \in N^\perp$, $\phi \in \H_N$ and $g' \in G$,
$(\p_N(s(n))\phi)(g')=\Xi_n(g')\phi(g')$.
\end{prop}

\begin{proof}
\begin{equation*}
\begin{aligned}
(\p_N(s(n))\phi)(g')
=&\; \phi(g's(n))\\
=&\; \phi(g's(n)g'^{-1}g')\\
=&\; \tau(g's(n)g'^{-1})\phi(g')\\
=&\; \tau(g's(n)g'^{-1}s(n)^{-1}s(n))\phi(g')\\
=&\; e(\pi(g'),n)\tau(s(n))\phi(g')\\
=&\; \Xi_n(g')\phi(g').
\end{aligned}
\end{equation*}
\end{proof}

\begin{defn}
If $F \in L^1(G//\pi^{-1}(N^\perp), \B(\H_\tau))$, and $n \in N^\perp$, we define
the $n$-th {\em Fourier-Wigner coefficient} of $F$ by
\begin{equation*}
(\alpha(F))(n)=\int_{\Omega_{N^\perp}} \<F(s(w)), \Xi_n(s(w))\> \, dw.
\end{equation*}
\end{defn}

\begin{prop}\label{P:FWC-uniqueness}
If $F \in L^1(G//\pi^{-1}(N^\perp), \B(\H_\tau))$ and $(\alpha(F))(n)=0$
for all $n \in N^\perp$, then $F=0$.
\end{prop}

\begin{proof}
Observe that if $m \in N$ and $n' \in S$, then
\begin{equation*}
\begin{aligned}
\tau(s(m+n'))
=&\; \tau(\psi(m,n')^{-1}s(m)s(n'))\\
=&\; \zeta(\psi(m,n')^{-1}s(m))s(n').
\end{aligned}
\end{equation*}
Therefore, for all $m \in N$ and $n' \in S$,
\begin{equation*}
\begin{aligned}
0
=&\; (\alpha(F))(m+n')\\
=&\; \int_{\Omega_{N^\perp}} \<F(s(w)), \Xi_{m+n'}(s(w))\> \, dw\\
=&\; \int_{\Omega_{N^\perp}} \<F(s(w)), \tau(s(m+n'))\> \bar{e(w, m+n')} \, dw\\
=&\; \bar{\zeta(\psi(m,n')^{-1}s(m))}
     \int_{\Omega_{N^\perp}} \<F(s(w)), \tau(s(n'))\>
                      \bar{e(w,n')} \, \bar{e(w,m)} \, dw.
\end{aligned}
\end{equation*}
Since the linear span of the set $\{e(\cdot,m) \st m \in N\}$ is weakly
dense in $L^\infty(\Omega_{N^\perp})$, it follows that for all $n' \in S$,
\begin{equation*}
\<F(s(w)), \tau(s(n'))\> \bar{e(w,n')} = 0
\end{equation*}
for almost all $w \in \Omega_{N^\perp}$.
Therefore, for all $n' \in S$,
\begin{equation*}
\<F(s(w)), \tau(s(n'))\> = 0
\end{equation*}
for almost all $w \in \Omega_{N^\perp}$.
Therefore, by Corollary \ref{C:finite-onb}
\begin{equation*}
F(s(w)) = 0
\end{equation*}
for almost all $w \in \Omega_{N^\perp}$.
\end{proof}

\begin{prop}\label{P:FWC-parseval}
The set $\{\Xi_n\}_{n\in N^\perp}$ is a complete orthonormal set in
$L^2(G//\pi^{-1}(N^\perp), \B(\H_\tau))$.
\end{prop}

\begin{proof}
By Corollary \ref{C:finite-inner-products},
\begin{equation*}
\begin{aligned}
\<\Xi_n, \Xi_{n'}\>
=&\; \int_{\Omega_{N^\perp}} \<\Xi_n(s(w)),\Xi_{n'}(s(w))\> \, dw\\
=&\; \int_{\Omega_{N^\perp}} e(w, n-n') \, dw \<\tau(s(n),\tau(s(n')\>\\
=&\;
\begin{cases}
a \zeta(\psi(n,-n') s(n-n')) \int_{\Omega_{N^\perp}} e(w,n-n') \, dw
& \text{if $n-n' \in N$},\\
0 & \text{if $n-n' \notin N$}
\end{cases}\\
=&\;
\begin{cases}
1 & \text{if $n=n'$},\\
0 & \text{if $n \ne n'$}
\end{cases}\\
\end{aligned}
\end{equation*}
Completeness follows from Proposition \ref{P:FWC-uniqueness}.
\end{proof}

Let $\B(\H_N)$ denote the set of all bounded linear operators on $\H_N$.
Suppose $Y=\phi\otimes\bar{\psi} \in \B(\H_N)$ is a rank-one operator.
Define $F_Y:G \to \B(\H_\tau)$ by $F_Y(g)=\phi(g)\otimes\bar{\psi(g)}$.
Then $F_Y$ is measurable, and for $h \in \pi^{-1}(N^\perp)$,
\begin{equation*}
\begin{aligned}
F_Y(hg)
=&\; \phi(hg)\otimes\bar{\psi(hg)}\\
=&\; \tau(h)\phi(g) \otimes \bar{\tau(h)\psi(g)}\\
=&\; \tau(h)(\phi(g) \otimes \bar{\psi(g)})\tau(h)^{-1},
\end{aligned}
\end{equation*}
so $F_Y \in L^0(G//\pi^{-1}(N^\perp), \B(\H_\tau))$.  Moreover,
\begin{equation*}
\begin{aligned}
\norm{F_Y}_1
=&\; \int_{\Omega_{N^\perp}}
     (\tr((\phi(s(w))\otimes\bar{\psi(s(w))})
         (\phi(s(w))\otimes\bar{\psi(s(w))})^*))^{1/2} \, dw \\
=&\; \int_{\Omega_{N^\perp}}
     (\tr((\phi(s(w))\otimes\bar{\psi(s(w))})
         (\psi(s(w))\otimes\bar{\phi(s(w))})))^{1/2} \, dw\\
=&\; \int_{\Omega_{N^\perp}}
     \norm{\psi(s(w))}
     (\tr(\phi(s(w))\otimes\bar{\phi(s(w))}))^{1/2} \, dw\\
=&\; \int_{\Omega_{N^\perp}}
     \norm{\psi(s(w))}\norm{\phi(s(w))} \, dw\\
\le &\; \norm{\psi}\norm{\phi} < \infty,
\end{aligned}
\end{equation*}
so $F_Y \in L^1(G//\pi^{-1}(N^\perp), \B(\H_\tau))$.
If $Y=\sum \phi_j \otimes \bar{\psi_j} \in \B(\H_N)$ is a finite
rank operator, then we define $F_Y=\sum F_{\phi_j \otimes \bar{\psi_j}}$,
and $F_Y \in L^1(G//\pi^{-1}(N^\perp), \B(\H_\tau))$.
The following obvious proposition will play a crucial role in the
proof of Theorem \ref{T:reformulation} in section \ref{S:benedicks}.

\begin{prop}\label{P:rank-comparison}
$\mathrm{rank}(F_Y(g)) \le \mathrm{rank}(Y)$ for almost all $g \in G$.
\end{prop}

If $X \in \B(\H)$, put $X_N=\Z X \Z^{-1}$, where $\Z:\H\to\H_N$
is a unitary intertwiner between $\p$ and $\p_N$.  Since $\Z$ is unique
up to a unitary constant, $X_N$ is determined uniquely by $X$.
Also $X_N \in \B(\H_N)$.  The operator theoretic properties of
$X$ are inherited by $X_N$, e.g. $X_N$ is of trace class
iff $X$ is of trace class, $X_N$ is injective iff $X$ is injective,
and $X_N$ is of finite rank iff $X$ is of finite rank, and moreover,
$\mathrm{rank}(X_N)=\mathrm{rank}(X)$.

Observe that if $X$ is a trace class operator on $\H$, then
\begin{equation*}
\tr(X_N \p_N(s(w))^* ) = \alpha(X)(w).
\end{equation*}

\begin{thm}\label{T:poisson}
If $X$ is a finite rank operator on $\H$, and $n \in N^\perp$, then
$\alpha(F_{X_N})(n)=\alpha(X)(n)$.
\end{thm}

\begin{proof}
It suffices to prove this for a rank one operator $X$.  In this case
$X_N=\phi\otimes\bar{\psi}$ for some $\phi,\psi \in \H_N$.  Therefore,
by Proposition \ref{P:trig-ops},
\begin{equation*}
\begin{aligned}
\alpha(X)(n)
=&\; \tr(X_N\p_N(s(n))^*)\\
=&\; \tr((\phi\otimes\bar{\psi})\p_N(s(n))^*)\\
=&\; \tr(\phi\otimes\bar{\p_N(s(n))\psi})\\
=&\; \<\phi,\p_N(s(n))\psi\>\\
=&\; \int_{\Omega_{N^\perp}}\<\phi(s(w)), (\p_N(s(n))\psi)(s(w))\> \, dw\\
=&\; \int_{\Omega_{N^\perp}}\<\phi(s(w)), \Xi_n(s(w))\psi(s(w))\> \, dw\\
=&\; \int_{\Omega_{N^\perp}}\tr(\phi(s(w))\otimes\bar{(\Xi_n(s(w))\psi(s(w)))}) \, dw\\
=&\; \int_{\Omega_{N^\perp}}\tr((\phi(s(w))\otimes\bar{\psi(s(w))})\Xi_n(s(w))^*) \, dw\\
=&\; \int_{\Omega_{N^\perp}}\tr(F_{X_N}(s(w))\Xi_n(s(w))^*) \, dw\\
=&\; \alpha(F_{X_N})(n).
\end{aligned}
\end{equation*}
\end{proof}

\section{The Benedicks argument}\label{S:benedicks}

We will now prove Theorem \ref{T:reformulation}.
Assume $X$ is a finite rank operator on $\H$.  Let
$B=\{w \in \RR^2 \st \alpha(X)(w) \ne 0\}$ and assume that $B$ has finite
measure.  Choose an isotropic lattice $N \subseteq \RR^2$ such that
a=$\mathrm{area}(\RR^2/N) > \mathrm{rank}(X)$.  Then
$N_v=(B-v) \cap N^\perp$ is a finite set for almost every $v \in \RR^2$.
For such $v$, let $X^v=X\p(s(v))^*$.  By Proposition \ref{P:rank-comparison},
$\mathrm{rank}(F_{X^v_N}(g')) \le \mathrm{rank}(X^v_N) =
 \mathrm{rank}(X^v)=\mathrm{rank}(X) < a$
for almost every $g' \in G$.  Since
\begin{equation*}
\begin{aligned}
\alpha(X^v)(w)
=&\; \tr(X^v\p(s(w))^*)\\
=&\; \tr(X\p(s(v))^*\p(s(w))^*)\\
=&\; \tr(X\p(s(w)s(v))^*)\\
=&\; \bar{\psi(w,v)}\tr(X\p(s(w+v))^*)\\
=&\; \bar{\psi(w,v)}\alpha(X)(w+v),
\end{aligned}
\end{equation*}
$\alpha(X^v)(w)=0$ if $w \notin B-v$, and so by Theorem \ref{T:poisson},
$\alpha(F_{X^v_N})(n)=\alpha(X^v)(n)=0$ if $n \in N \setminus N_v$.
Therefore, by Proposition \ref{P:FWC-uniqueness} and
Proposition \ref{P:FWC-parseval},
\begin{equation*}
F_{X^v_N} = \sum_{n\in N_v} \alpha(X^v)(n) \Xi_n.
\end{equation*}

According to \cite[Theorem 1.2]{Linnel},
$\sum_{n \in N_v} \alpha(X^v)(n) \p(s(n))$ is injective,
unless it is zero.  Therefore 
\begin{equation*}
\sum_{n \in N_v} \alpha(X^v)(n) \p_N(s(n))
\end{equation*}
is injective, unless it is zero.
However, by Proposition \ref{P:trig-ops}, for all
$\phi \in \H_N$ and $g' \in G$,
\begin{equation*}
F_{X^v_N}(g')\phi(g') =
\sum_{n\in N_v} \alpha(X^v)(n) \Xi_n(g')\phi(g') =
\sum_{n\in N_v} \alpha(X^v)(n)(\p_N(s(n))\phi)(g').
\end{equation*}
Therefore
$F_{X^v_N}(g')$ has rank $a$, for almost every $g' \in G$, unless
it is identically zero.  So in fact $F_{X^v_N} \equiv 0$, and so
$\alpha(X^v)(n)=0$ for all $n \in N$.  Therefore
$\alpha(X)(n+v)=0$ for all $n \in N$ and almost every $v \in \RR^2$.
Therefore $\alpha(X)=0$ almost everywhere, and so $X=0$.

\bibliographystyle{amsplain}
\bibliography{v11-btwt}

\end{document}